\documentclass[11pt]{amsart}

\usepackage{amsmath}
\usepackage{amssymb}
\usepackage{mathrsfs}

\newtheorem{theorem}{Theorem}[section]
\newtheorem{claim}[theorem]{Claim}

\newtheorem{proposition}[theorem]{Proposition}
\newtheorem{corollary}[theorem]{Corollary}

\theoremstyle{definition}
\newtheorem{definition}[theorem]{Definition}

\newtheorem{question}[theorem]{Question}

\theoremstyle{remark}

\newcount\skewfactor
\def\mathunderaccent#1#2 {\let\theaccent#1\skewfactor#2
\mathpalette\putaccentunder}
\def\putaccentunder#1#2{\oalign{$#1#2$\crcr\hidewidth
\vbox to.2ex{\hbox{$#1\skew\skewfactor\theaccent{}$}\vss}\hidewidth}}
\def\name{\mathunderaccent\tilde-3 }

\def\smallbox#1{\leavevmode\thinspace\hbox{\vrule\vtop{\vbox
   {\hrule\kern1pt\hbox{\vphantom{\tt/}\thinspace{\tt#1}\thinspace}}
   \kern1pt\hrule}\vrule}\thinspace}

\def\l{{\langle}}
\def\r{{\rangle}}

\newcommand{\cf}{{\rm cf}}

\def\qedref#1{$\qed_{\reforiginal{#1}}$}

\title{Kurepa trees and the failure of the Galvin property}

\author{Tom Benhamou}
\address{School of Mathematical Sciences, Raymond and Beverly Sackler Faculty of Exact Science, Tel-Aviv University, Ramat Aviv 69978, Israel}
\email{tombenhamou@tauex.tau.ac.il}

\author{Shimon Garti}
\address{Einstein Institute of Mathematics,
 The Hebrew University of Jerusalem,
 Jerusalem 91904, Israel}
\email{shimon.garty@mail.huji.ac.il}

\author{Saharon Shelah}
\address{Institute of Mathematics
 The Hebrew University of Jerusalem,
 Jerusalem 91904, Israel
 and  Department of Mathematics
 Rutgers University
 New Brunswick, NJ 08854, USA}
\email{shelah@math.huji.ac.il}
\urladdr{http://www.math.rutgers.edu/\char`\~shelah}
\thanks{The research was supported by Israel Science Fountdation Grant no. 1838/19. This is publication 1219 of the third author}

\begin{document}
\let\labeloriginal\label
\let\reforiginal\ref
\def\ref#1{\reforiginal{#1}}
\def\label#1{\labeloriginal{#1}}

\begin{abstract}
We force the existence of a non-trivial $\kappa$-complete ultrafilter over $\kappa$ which fails to satisfy the Galvin property. This answers a question asked in \cite{bg}.
\end{abstract}

\maketitle

\newpage

\section{Introduction}

Let $\kappa$ be strongly regular (that is, $\kappa=\kappa^{<\kappa}$) and uncountable.
Let $\mathscr{F}$ be a normal filter over $\kappa$.
A delightful theorem of Galvin says that if one considers a family $\mathcal{C}=\{C_\alpha:\alpha\in\kappa^+\}\subseteq\mathscr{F}$ then one can always find a subfamily $\{C_{\alpha_i}:i\in\kappa\}$ of $\mathcal{C}$ such that $\bigcap\{C_{\alpha_i}:i\in\kappa\}\in\mathscr{F}$.
The statement and the proof were published in \cite{MR0369081}.

One may wonder whether the assumption $\kappa=\kappa^{<\kappa}$ is necessary.
The answer is yes, and the failure of the Galvin property is consistent as proved by Abraham and Shelah in \cite{MR830084}.
The failure of the Galvin property was further investigated in \cite{MR3604115}, \cite{MR3787522} and \cite{bgp}.
Of course, $\kappa<\kappa^{<\kappa}$ in any such model, hence if $\kappa$ is strongly inaccessible then the Galvin property holds at every normal filter over $\kappa$.
Following the notation of \cite{bgp} we denote by ${\rm Gal}(\mathscr{F},\kappa,\kappa^+)$ the following statement:
$$\forall \{A_\alpha:\alpha\in\kappa^+\}\subseteq\mathscr{F}, \ \exists I\in[\kappa^+]^{\kappa}, \ \cap_{\alpha\in I}A_\alpha\in \mathscr{F}.$$

Galvin's theorem applies to any normal filter.
In particular, if $\mathscr{U}$ is a normal \emph{ultrafilter} over a measurable cardinal $\kappa$ then ${\rm Gal}(\mathscr{U},\kappa,\kappa^+)$ holds true.

The normality of the filter plays an important r\^ole in Galvin's proof.
Essentially, given a family $\mathcal{C}=\{C_\alpha:\alpha\in\kappa^+\}$ one can isolate a subfamily $\{C_{\alpha_i}:i\in\kappa\}$ of $\mathcal{C}$ whose intersection equals its diagonal intersection up to some $\mathscr{F}$-negligible set.
Normality, therefore, is crucial here since it makes sure that the filter is closed under diagonal intersections.

However, the scope of Galvin's property is wider.
For example, the first author and Gitik proved in \cite{bg} that if $\mathscr{U}$ is a $p$-point or a product of $p$-points over $\kappa$ then ${\rm Gal}(\mathscr{U},\kappa,\kappa^+)$ holds.
Following this statement they posed several questions about the applicability of Galvin's theorem to non-normal ultrafilters.
One of the interesting cases is a measurable cardinal $\kappa$, since such a cardinal carries many $\kappa$-complete ultrafilters which are not necessarily normal.

\begin{question}
\label{qmq} Let $\kappa$ be a measurable cardinal.
\begin{enumerate}
\item [$(\aleph)$] Is it consistent that $\neg{\rm Gal}(\mathscr{U},\kappa,\kappa^+)$ for some $\kappa$-complete ultrafilter over $\kappa$?
\item [$(\beth)$] Is it consistent that $\neg{\rm Gal}(\mathscr{U},\kappa,\kappa^+)$ where $\mathscr{U}$ is $\kappa$-complete and extends the club filter of $\kappa$?
\end{enumerate}
\end{question}

It has been shown in \cite{bg} that consistently $\kappa$ is a measurable cardinal and ${\rm Gal}(\mathscr{U},\kappa,\kappa^+)$ holds at every $\kappa$-complete ultrafilter over $\kappa$.
To see this, consider Solovay's inner model $L[\mathscr{U}]$ where $\mathscr{U}$ is a normal ultrafilter over $\kappa$.
Kunen proved in \cite{MR277346} that if $\mathscr{W}$ is a $\kappa$-complete ultrafilter over $\kappa$ in $L[\mathscr{U}]$ then $\mathscr{W}$ is Rudin-Keisler equivalent to $\mathscr{U}^n$ for some $n\in\omega$.
Now if $\mathscr{U}$ and $\mathscr{V}$ are Rudin-Keisler equivalent then ${\rm Gal}(\mathscr{U},\kappa,\kappa^+)$ holds iff ${\rm Gal}(\mathscr{V},\kappa,\kappa^+)$ holds.
Since ${\rm Gal}(\mathscr{U}^n,\kappa,\kappa^+)$ is true for every $n\in\omega$, the above statement follows.

The main result of this paper is the consistency of the opposite situation.
Namely, it is consistent that $\mathscr{U}$ is a $\kappa$-complete ultrafilter over $\kappa$ and ${\rm Gal}(\mathscr{U},\kappa,\kappa^+)$ fails.
Moreover, it is possible to force this failure for some $\mathscr{U}$ which extends the club filter of $\kappa$.
This result shows that normality and $\kappa$-completeness differ with respect to basic combinatorial properties, as shown e.g. in \cite{MR3902803}.

Our notation is standard for the most part.
If $\kappa=\cf(\kappa)<\mu$ then $S^\mu_\kappa=\{\delta\in\mu:\cf(\delta)=\kappa\}$.
If $\mu\geq\cf(\mu)>\omega$ then this set is stationary in $\mu$.
We employ the Jerusalem forcing notation, so $p\leq_{\mathbb{P}}q$ means that $p$ is weaker than $q$.
If $\kappa=\cf(\kappa)>\aleph_0$ then the club filter over $\kappa$ is denoted by $\mathscr{D}_\kappa$.
If $\kappa$ is supercompact then a Laver diamond for $\kappa$ is a function $h:\kappa\rightarrow{V_\kappa}$ which enjoys the following property.
For every set $x$ and every sufficiently large $\lambda>\kappa$ there is a fine and normal meausre $\mathscr{U}$ over $[\lambda]^{<\kappa}$ such that $\jmath_{\mathscr{U}}(h)(\kappa)=x$.
Finally, the polarized relation $\binom{\alpha}{\beta}\rightarrow\binom{\gamma}{\delta}_\theta$ means that for every $c:\alpha\times\beta\rightarrow\theta$ one can find $A\in[\alpha]^\gamma, B\in[\beta]^\delta$ and $i\in\theta$ such that $c\upharpoonright(A\times B)$ is constantly $i$.
According to Hajnal in \cite{MR281605}, this notation is due to Galvin.
For basic background in polarized partition relations (which play a central role in our proof) we suggest \cite{MR3075383}.

\newpage

\section{Preliminaries}

In this section we collect some information about Kurepa trees.
A Kurepa tree on $\aleph_1$ is a tree $\mathscr{T}$ of height $\aleph_1$ such that every level of $\mathscr{T}$ is countable and there are at least $\aleph_2$-many cofinal branches of $\mathscr{T}$.
As depicted in \cite[Chapter 31]{ezekiel}, it exceeded in stature all the trees of the field, its branches multiplied and its boughs grew long.

We are interested in the generalization of this concept to strongly inaccessible cardinals.
If one imitates the classical definition of Kurepa trees then the tree ${}^{<\kappa}2$ is Kurepa whenever $\kappa$ is strongly inaccessible, thus the concept of Kurepa trees becomes uninteresting.
We shall use the traditional substitute:

\begin{definition}[Slim Kurepa trees]
\label{defslim}
A $\kappa$-tree $\mathscr{T}$ is a slim Kurepa tree iff $|\mathcal{L}_\beta(\mathscr{T})|\leq|\beta|$ for every $\beta\in\kappa$ and the cardinality of the set of cofinal branches of $\mathscr{T}$ is at least $\kappa^+$.
\end{definition}

Before proceeding let us exclude from the discussion a seemingly stronger concept.
Call $\mathscr{T}$ \emph{very slim} if $|\mathcal{L}_\beta(\mathscr{T})|<|\beta|$ for every $\beta\in\kappa$.
Let us show that the number of cofinal branches in such trees must be small.

\begin{claim}
\label{clmveryslim} Suppose that $\kappa=\cf(\kappa)>\aleph_0$.
Then there are no very slim $\kappa$-Kurepa trees.
\end{claim}

\par\noindent\emph{Proof}. \newline
Let $\mathscr{T}$ be a very slim $\kappa$-tree.
We intend to prove that the cardinality of the set of cofinal branches of $\mathscr{T}$ is strictly less than $\kappa$.
For every $\beta\in\kappa$ let $\theta_\beta=|\mathcal{L}_\beta(\mathscr{T})|$, so $\theta_\beta<|\beta|$.
By Fodor's lemma there are a stationary $S\subseteq\kappa$ and a fixed cardinal $\theta\in\kappa$ such that $\beta\in S\Rightarrow\theta_\beta=\theta$.
Choose such a pair $(S,\theta)$ where $\theta$ is minimal.
Let $\mathcal{B}$ be the set of cofinal branches of $\mathscr{T}$.
We claim that $|\mathcal{B}|\leq\theta$.

Assume towards contradiction that $|\mathcal{B}|>\theta$ and let $\{b_\alpha:\alpha\in\theta^+\}\subseteq\mathcal{B}$ be a set of pairwise distinct cofinal branches.
For every pair of ordinals $\{\alpha,\delta\}\in[\theta^+]^2$ let $\beta_{\alpha\delta}\in\kappa$ be so that $b_\alpha(\beta_{\alpha\delta})\neq b_\delta(\beta_{\alpha\delta})$.
Let $\xi=\bigcup\{\beta_{\alpha\delta}:\{\alpha,\delta\}\in[\theta]^2\}$.
Notice that $\xi\in\kappa$ since $\kappa=\cf(\kappa)>\theta^+$.
Choose $\beta\in S$ such that $\beta>\xi$.
For every $\alpha\in\theta^+$ let $x_\alpha=b_\alpha\cap\mathcal{L}_\beta(\mathscr{T})$.
By the above choices if $\alpha\neq\delta$ then $x_\alpha\neq x_\delta$.
However, $|\mathcal{L}_\beta(\mathscr{T})|=\theta<\theta^+$, so this is impossible.

\hfill \qedref{clmveryslim}

Back to the concept of slim Kurepa trees, one can force such trees at smallish large cardinals.
In what follows, we always assume that a tree is downward-closed.
Let $\mathbb{K}$ be the following forcing notion described in \cite{MR2768691}, which adds a slim Kurepa tree.
A condition $p=(t^p,f^p)$ in $\mathbb{K}$ is a pair in which $t^p$ is a normal tree of height $\beta+1$ (where $\beta\in\kappa$) and $|\mathcal{L}_\alpha(t^p)|\leq|\alpha|$ whenever $\alpha\leq\beta$, so $t^p$ approximates a slim $\kappa$-tree.
In addition, $f^p$ is a partial function from $\kappa^+$ into $\mathcal{L}_\beta(t^p)$ (where $\beta$ is the index of the maximal level of $t^p$) and $|{\rm dom}(f^p)|\leq|\beta|$.
If $p,q\in\mathbb{K}$ then $p\leq_{\mathbb{K}}q$ iff $t^p$ is a subtree of $t^q, {\rm dom}(f^p)\subseteq{\rm dom}(f^q)$ and $f^p(\delta)\leq_{t^q}f^q(\delta)$ whenever $\delta\in{\rm dom}(f^p)$.

It is easy to see that if $\kappa^{<\kappa}=\kappa$ then $\mathbb{K}$ is $\kappa$-complete and $\kappa^+$-cc.
Hence if $G\subseteq\mathbb{K}$ is $V$-generic then no cardinals are collapsed in $V[G]$.
Let $\mathscr{T}=\bigcup\{t^p:p\in{G}\}$.
By density arguments one can see that $\mathscr{T}$ is a slim $\kappa$-Kurepa tree, and if $\kappa$ is strongly inaccessible then it remains so in $V[G]$.
Moreover, one can force such trees and preserve the property of weak compactness.

However, if $\kappa$ is sufficiently large then no such trees can live over $\kappa$.
If $\kappa$ is ineffable then there are no slim $\kappa$-Kurepa trees.
Actually, in $L$ this fact characterizes ineffability, see \cite{MR750828}.
A fortiori, if $\kappa$ is measurable then every $\kappa$-tree with $\kappa^+$ cofinal branches is not slim.
To see this, suppose that $\mathscr{T}$ is such a tree and let $\jmath:V\rightarrow{M}$ be an elementary embedding with ${\rm crit}(\jmath)=\kappa$.
By elementarity, $\jmath(\mathscr{T})$ is a $\jmath(\kappa)$-tree in $M$.
But $\jmath(\mathscr{T})$ is not slim since $|\mathcal{L}_\kappa(\jmath(\mathscr{T}))|\geq\kappa^+$, a fact which follows from the existence of $\kappa^+$-many cofinal branches of $\mathscr{T}$ in $V$.
By elementarity, again, $\mathscr{T}$ is not slim in $V$.

This simple observation sets down a constraint on slim Kurepa trees at measurable cardinals, but it also opens a window to a possible weakening of this concept which can live happily with measurability and even stronger large cardinal properties.
The point is that the non-slimness seems to concentrate on certain levels below $\jmath(\kappa)$, so by restricting the slimness to non-problematic levels we may hope to fix this issue.

\begin{definition}
Let $\kappa$ be a strongly inaccessible cardinal and let $S$ be a stationary subset of $\kappa$.
A slim $S$-Kurepa tree is a tree $\mathscr{T}\subseteq {}^{<\kappa}\kappa$ such that for every $\alpha\in S$, $|\mathcal{L}_\alpha(S)|\leq |\alpha|$ and $\sup(\mathscr{T}):=\{t\in{}^{\kappa}\kappa:\forall\alpha<\kappa, t\upharpoonright\alpha\in\mathscr{T}\}$ is a set of size at least $\kappa^+$. Also $\mathscr{T}$ is called a  stationary-slim Kurepa tree if for some stationary set $S\subseteq \kappa$, $\mathscr{T}$ is a slim $S$-Kurepa tree.
\end{definition}

Let us briefly describe the idea behind our next definition.
Adding a Cohen subset to a measurable cardinal $\kappa$ may destroy measurability.
But adding it on top of an Easton product which adds a Cohen subset at every strongly inaccessible $\alpha\in\kappa$ will preserve measurability.
This can be viewed as a baby-version of Laver's indestructibility from \cite{MR0472529}.
While Laver wanted to take care of every $\kappa$-directed-closed forcing notion, it is possible to focus on a narrower class of forcing notions.
In the following we define a forcing notion which applies this idea to the poset which adds a slim Kurepa tree.

\begin{definition}[The forcing for adding a stationary-slim Kurepa tree]
Let $\kappa$ be a regular cardinal and $S$ a stationary subset of $\kappa$.
The forcing $\mathbb{Q}(S)$ is the forcing for adding a stationary subset of $S$, namely, $\mathbb{Q}=\{X\subseteq S:|X|<\kappa\}$ where the order is end extension which is denoted by $\leq_{end}$.
Define the forcing $\mathbb{K}(S)$ to be the forcing notion which consists of triples $\l X,t,f\r$ where:
\begin{enumerate}
\item [$(a)$] $X\in \mathbb{Q}(S)$.
\item [$(b)$] $t$ is a normal tree of height $\beta+1<\kappa$, $\beta+1\leq\sup(X)$.
\item [$(c)$] For every $\alpha\in X\cap \beta+1$, $|\mathcal{L}_\alpha(t)|\leq |\alpha|$.
\item [$(d)$] $f:\kappa^+\rightarrow \mathcal{L}_\beta(t)$ is a partial function, and $|f|\leq |\beta|$.
\end{enumerate}
For the order, we set $\l X,t,f\r\leq \l Y,s,g\r$ iff:
\begin{enumerate}
\item [$(\alpha)$] $X\leq_{end} Y$.
\item [$(\beta)$] $s\restriction \beta+1=t$.
\item [$(\gamma)$] $Dom(f)\subseteq Dom(g)$.
\item [$(\delta)$] For every $\alpha\in Dom(f)$, $f(\alpha)\leq_s g(\alpha)$.
\end{enumerate}
Denote the projection to each coordinate by $\pi_1(\l a,b,c\r)=a,\ \pi_2(\l a,b,c\r)=b,\ \pi_3(\l a,b,c\r)=c$.
\end{definition}

The following is clear:

\begin{proposition}
\label{propchaincondition}
If $\kappa^{<\kappa}=\kappa$ then $\langle\mathbb{K}(S),\leq\rangle$ is $\kappa$-closed and $\kappa^+$-cc.
\end{proposition}

The proposition below shows that $\mathbb{K}(S)$ adds a stationary set and a Kurepa tree which is slim in this set.

\begin{proposition}\label{propshooting}
Let $S\subseteq S^\kappa_\omega$ be stationary and let $G$ be $V$-generic for $\mathbb{K}(S)$, define $S^*=\bigcup\{\pi_1(p):p\in G\}$ and $\mathscr{T}=\bigcup\{\pi_2(p)\mid p\in G\}$.
Then $S^*$ is stationary at $\kappa$ and $\mathscr{T}$ is a slim $S^*$-Kurepa tree.
\end{proposition}

\par\noindent\emph{Proof}. \newline
Let $C\subseteq\kappa$ be any club in $V[G]$, and let $\name{C}$ be a $\mathbb{K}(S)$-name for it.
Choose $p\in G$ which forces that $\name{C}$ is a club at $\kappa$.
We proceed by a density argument.
Let $q_0\geq p$ be any condition. Let us define inductively sequences $(q_\beta : \beta<\kappa),(\delta_\beta:\beta<\kappa)$ such that:\begin{enumerate}
    \item $(q_{\beta} : \beta<\kappa)$ is increasing and continuous, namely $\pi_1(q_\beta)=\cup_{j<\beta}\pi_1(q_j)$.
    \item $\sup(\pi_1(q_{\beta}))\leq \delta_\beta$ and $q_{\beta+1}\Vdash \delta_{\beta}\in\name{C}$.
    \item $(\delta_\beta:\beta<\kappa)$ is strictly increasing and continuous.
\end{enumerate}
The condition $q_0$ was already defined. Since $\mathbb{K}(S)$ is $\kappa$-closed, the set $E_{q_0}:=\{\alpha\in\kappa: \exists q'\geq q_0, q'\Vdash \alpha\in\name{C}\}$ contains a club subset of $\kappa$, hence for some $\delta_0\in E_{q_0}-\sup(\pi_1(q_0))$ and $q_0\leq q_1$, we have $q_1\Vdash \delta_0\in \name{C}$.
Assume we have defined $q_\alpha,\delta_\alpha$ for
$\alpha<\beta<\kappa$, satisfying $(1)-(3)$. For limit $\beta$, let
$q_\beta$ be a least upper bound of $(q_\alpha:\alpha<\beta)$, in particular $\pi_1(q_{\beta})=\cup_{\alpha<\beta}\pi_1(q_{\alpha})$, hence by induction  $\delta_\beta:=\sup(\delta_\alpha:\alpha<\beta)\geq \sup(\pi_1(q_{\beta}))$. Note that since $q_{\beta}\geq q_\alpha$ for each $\alpha<\beta$, $q_\beta\Vdash \delta_\alpha\in\name{C}$ and since $q_{\beta}\Vdash \name{C}$ is a club, we have $q_{\beta}\Vdash \delta_\beta\in\name{C}$.
At successor step $\beta+1$, define $E_{q_\beta}$ as before, but using the condition $q_{\beta}$ instead of $q_0$, and find $\delta_{\beta+1}\in  E_{q_\beta}-\sup(\pi_1(q_{\beta})\cup\{\delta_{\beta}+1\})$ and $q_{\beta+1}\geq q_\beta$ such that $q_{\beta+1}\Vdash \delta_{\beta+1}\in \name{C}$. This concludes the construction.

By $(3)$, the set $C^*=\{\delta_\beta: \beta\text{ is limit}, \beta<\kappa\}$ forms a club.
By the stationarity of $S$, there is some limit $\beta$ such that $\delta_\beta\in S\cap C^*$.
Consider $q_\beta$, then by $(2)$, $\sup(\pi_1(q))\leq \delta_\beta$. Extend $q$ to $q^+$ such that $\delta_\beta\in \pi_1(q^+)$, this is possible since $\delta_{\beta}\in S$ and above $\sup(\pi_1(q^+)$. Thus $q^+\Vdash \delta_{\beta}\in \name{S}^*$. As in the limit step, note that $q^+\geq q_{\beta}\geq q_{i+1}$ for every $i<\beta$, hence by $(2)$ $q^+\Vdash \delta_i\in\name{C}$, hence $q^+\Vdash \delta_{\beta}\in \name{C}$.
We conclude that $q^+\Vdash \delta_\beta\in \name{C}\cap\name{S}^*$.
By density this means that $ C\cap S^*\neq\varnothing$.

For the second part, clearly $\mathscr{T}$ is a normal tree such that for each $\alpha\in S^*$, $\mathcal{L}_\alpha(\mathscr{T})$ is of size at most $\alpha$.
Let us argue that there are $\kappa^+$-many branches in $\mathscr{T}$.
Indeed for every $\alpha\in\kappa^+$ the set $\{f(\alpha):\exists t,X, \l X,t,f\r\in G\wedge \alpha\in Dom(f)\}$ is a branch in $\mathscr{T}$.
By density we can make these branches distinct, so we are done.

\hfill \qedref{propshooting}

\newpage

\section{Kurepa meets Galvin}

In this section we prove the main result of the paper.
Our strategy is to prove two incompatible statements.
One of them is a negative combinatorial claim which follows from the existence of a certain slim $\kappa$-Kurepa tree.
The opposite one is a positive statement which follows from the assumption that every $\kappa$-complete ultrafilter over $\kappa$ satisfies the Galvin property.
The pertinent combinatorial property is a generalization of the following (unpublished) result also due to Galvin.

\begin{claim}
\label{conclmgalvin}
If there exists a Kurepa tree then $\binom{\omega_2}{\omega_1}\nrightarrow\binom{2}{\omega_1}_\omega$.
\end{claim}

\par\noindent\emph{Proof}. \newline
Suppose that $f,g:\omega_1\rightarrow\omega$.
We shall say that $f$ and $g$ are \emph{almost disjoint} iff $|\{\beta\in\omega_1:f(\beta)=g(\beta)\}|\leq\aleph_0$.
Let $\mathscr{T}$ be a Kurepa tree.
We intend to build a collection $\mathscr{F} = \{f_\alpha:\alpha\in\omega_2\}\subseteq{}^{\omega_1}\omega$ such that $\mathscr{F}$ is a family of pairwise almost disjoint functions.

For every $\beta\in\omega_1$ fix an enumeration $(t_{\beta n}:n\in\omega)$ of the elements of $\mathcal{L}_\beta(\mathscr{T})$.
Let $(b_\alpha:\alpha\in\omega_2)$ be an enumeration of the $\omega_1$-branches of $\mathscr{T}$ (or some of them, if there are more than $\aleph_2$ branches).
For every $\alpha\in\omega_2$ define $f_\alpha:\omega_1\rightarrow\omega$ as follows:
$$
f_\alpha(\beta)=m \Leftrightarrow b_\alpha\cap\mathcal{L}_\beta(\mathscr{T})=t_{\beta m}.
$$
Notice that if $\alpha_0<\alpha_1<\omega_2$ then for some $\beta_0\in\omega_1$ we have $b_{\alpha_0}\upharpoonright\beta_0 = b_{\alpha_1}\upharpoonright\beta_0$ and $b_{\alpha_0}(\beta)\neq b_{\alpha_1}(\beta)$ whenever $\beta\in[\beta_0,\omega_1)$.
By the definition of our functions we see that $f_{\alpha_0}(\beta)\neq f_{\alpha_1}(\beta)$ for every $\beta\in[\beta_0,\omega_1)$.

Letting $\mathscr{F}=\{f_\alpha:\alpha\in\omega_2\}$ we may conclude that $\mathscr{F}$ is almost disjoint.
All we need now is to convert such a family to a coloring which exemplifies the negative relation $\binom{\omega_2}{\omega_1}\nrightarrow\binom{2}{\omega_1}_\omega$.
To this end, define a coloring $c:\omega_2\times\omega_1\rightarrow\omega$ as follows:
$$
c(\alpha,\beta)=f_\alpha(\beta).
$$
If $\alpha_0<\alpha_1<\omega_2$ and $B\in[\omega_1]^{\omega_1}$ then for a sufficiently large $\beta\in B$ we have $c(\alpha_0,\beta)=f_{\alpha_0}(\beta)\neq f_{\alpha_1}(\beta)=c(\alpha_1,\beta)$, so we are done.

\hfill \qedref{conclmgalvin}

The above claim generalizes, verbatim, to every pair of successor and double successor cardinals.
Thus, if there exists a $\kappa^+$-Kurepa tree then $\binom{\kappa^{++}}{\kappa^+}\nrightarrow\binom{2}{\kappa^+}_\kappa$.
In order to apply a similar idea to a limit cardinal and its successor one has to slightly modify the statement.
Rather than a subset of the small component of full size, we require a stationary set.

\begin{proposition}
\label{propgenstat} Suppose that $S\subseteq S^{\kappa}_{\aleph_0}$ is stationary and that there is a slim $S$-Kurepa tree on $\kappa$.
Then there is a coloring $c:\kappa^+\times S\rightarrow \omega$ with no monochromatic product of the form $A\times S'$ for $A\in[\kappa^+]^\kappa$ and $S'\subseteq S$ stationary.
\end{proposition}

\par\noindent\emph{Proof}. \newline
Let $\mathscr{T}$ be a slim $S$-Kurepa tree on $\kappa$.
For every $\alpha\in S$, $|\mathcal{L}_\alpha(\mathscr{T})|\leq \alpha$, and since $\cf(\alpha)=\aleph_0$, we can  partition $\mathcal{L}_\alpha(\mathscr{T})=\biguplus_{n\in\omega}I^{\alpha}_n$ such that $|I^\alpha_n|=\theta^\alpha_n<\alpha$.
Let $\{b_i:i\in\kappa^+\}$ be any collection of $\kappa^+$-many branches of $\mathscr{T}$ (which exists by the assumption that $\mathscr{T}$ is Kurepa).

Define the coloring $c:\kappa^+\times S\rightarrow \omega$ by $c(i,\alpha)=n$ for the unique $n$ such that $b_i\upharpoonright\alpha\in I^\alpha_n$.
Let us argue that there is no monochromatic product of the alleged form under $c$.
Suppose otherwise that $A\in[\kappa^+]^\kappa$ and $S'\subseteq S$ is stationary and $c\upharpoonright(A\times S')$ is constantly $n^*$.
The map $\alpha\mapsto \theta^\alpha_{n^*}$ is regressive on $S'$, hence by Fodor's lemma, there is a stationary $T'\subseteq S'$ such that for every $\alpha,\beta\in T'$, $\theta^{\alpha}_{n^*}=\theta^{\beta}_{n^*}=\theta$.

Let $A'\subseteq A$ be any subset of size $\theta^{+}$ and for every $i,j\in A'$ let $l_{i,j}\in\kappa$ be a level such that $b_j\upharpoonright l_{i,j}\neq b_i\upharpoonright l_{i,j}$.
By the regularity of $\kappa$, $\delta=\sup_{i,j\in A'}l_{i,j}\in\kappa$.
Since $T'$ is stationary one can find $\alpha\in T'-\delta$.
On the one hand, we know that $|I^{\alpha}_{n^*}|=\theta$.
On the other hand $\{b_i\upharpoonright \alpha:i\in A'\}\subseteq I^{\alpha}_{n^*}$, a contradiction.

\hfill \qedref{propgenstat}

The following proposition shows that the violation of this combinatorial property implies that the Galvin property must fail for a suitable ultrafilter.

\begin{proposition}
\label{propimplication}
Suppose that $\mathscr{U}$ is a $\sigma$-complete ultrafilter over $\kappa$, and ${\rm Gal}(\mathscr{U},\kappa,\kappa^+)$ holds.
Then for every coloring $c:\kappa^+\times S\rightarrow\omega$ where $S\in \mathscr{U}$, there are $A\in[\kappa^+]^\kappa$ and $S'\subseteq S$, $S'\in \mathscr{U}$ such that $c\upharpoonright(A\times S')$ is constant.
\end{proposition}

\par\noindent\emph{Proof}. \newline
For every $\alpha\in\kappa^+$ and every $n\in\omega$ let $S^{\alpha}_n=\{\beta\in S:c(\alpha,\beta)=n\}$.
Since $S=\bigcup_{n\in\omega} S^{\alpha}_n$ and $\mathscr{U}$ is $\sigma$-complete, there is $n_\alpha\in\omega$ such that $S^{\alpha}_{n_\alpha}\in \mathscr{U}$.  By the pigeon-hole principle, find $X\subseteq\kappa^+$, $|X|=\kappa^+$ such that for every $\alpha,\beta\in X$, $n_\alpha=n_\beta=m$.
Apply ${\rm Gal}(\mathscr{U},\kappa,\kappa^+)$ to the sequence $\l S^{\alpha}_{m}:\alpha\in X\r$ to find $A\in [X]^\kappa$, such that $S'=\bigcap_{\alpha\in A}S^{\alpha}_{m}\in\mathscr{U}$.
It follows that $c\upharpoonright(A\times S')=m$, so the proof is accomplished.

\hfill \qedref{propimplication}

Let us summarize the above propositions:

\begin{corollary}
\label{cor1}
Suppose $S\subseteq S^\kappa_{\aleph_0}$ is stationary and there is a slim $S$-Kurepa tree.
Then $\neg{\rm Gal}(\mathscr{U},\kappa,\kappa^+)$ holds for every $\sigma$-complete ultrafilter $\mathscr{U}$ which extends $\mathscr{D}_\kappa\cup\{S\}$.
\end{corollary}

\hfill \qedref{cor1}

It remains to prove that there is a model in which one can find a stationary subset $S\subseteq S^{\kappa}_{\aleph_0}$, a slim $S$-Kurepa tree and a $\kappa$-complete ultrafilter which extends $\mathscr{D}_\kappa\cup\{S\}$.
In \cite[Page 1033]{MR2987148}, P. L\"{u}cke indicates that under large cardinal assumptions it is possible to force a weak-Kurepa tree and keep $\kappa$ supercompact.
The concept of a weak Kurepa tree is parallel to what we call a slim Kurepa tree.
Let us provide a proof for our version of slimness in stationary sets which consist of ordinals of countable cofinality.
We indicate that countable cofinality is not essential.
That is, one can force slim $S$-Kurepa trees while preserving large cardinal properties for many other stationary sets $S$.
Likewise, one can prove similar negative combinatorial statements from the existence of such a tree, though the number of colors depends, usually, on the nature of the stationary set.

\begin{theorem}
\label{thmslimsupercompact}
Let $\kappa$ be a supercompact cardinal.
There is a forcing notion $\mathbb{S}$ such that in $V^{\mathbb{S}}$ there is a stationary set $S\subseteq S^{\kappa}_{\aleph_0}$, there is a slim $S$-Kurepa tree and $\mathscr{D}_\kappa\cup\{S\}$ can be extended to a $\kappa$-complete ultrafilter.
\end{theorem}

\par\noindent\emph{Proof}. \newline
Let $h:\kappa\rightarrow V_\kappa$ be a Laver-diamond function, as defined in \cite{MR0472529}.
The forcing $\mathbb{S}$ is defined as the iteration $\mathbb{P}_\kappa\ast\name{\mathbb{R}}$ where $\mathbb{P}_\kappa$ is an Easton support iteration $\l \mathbb{P}_\beta,\name{\mathbb{Q}}_\alpha: \beta\leq\kappa, \alpha<\kappa\r$ such that for every $\alpha\in\kappa$,  $\name{\mathbb{Q}}_\alpha$ is trivial, unless $\alpha$ is inaccessible and $h(\alpha)$ is a $\mathbb{P}_\alpha$-name for the forcing $\mathbb{K}(S^{\alpha}_{\aleph_0})$ in which case $\name{\mathbb{Q}}_\alpha=h(\alpha)$.

The forcing $\mathbb{R}$ is simply $\mathbb{K}(S^{\kappa}_{\aleph_0})$.
Let $G_{\mathbb{P}_\kappa}\ast G_{\mathbb{R}}$ be a $V$-generic set for $\mathbb{P}_\kappa\ast \name{\mathbb{R}}$.
By proposition \ref{propshooting},
in $V[G_{\mathbb{P}_\kappa}\ast G_{\mathbb{R}}]$ there is a stationary $S\subseteq S^{\kappa}_{\aleph_0}$ and a slim $S$-Kurepa tree.
It remains to show that there is a $\kappa$-complete ultrafilter extending $\mathscr{D}_\kappa\cup\{S\}$.

Since $h$ is a Laver-diamond, there is an embedding embedding $j:V\rightarrow M$ such that $M^{2^\kappa}\subseteq M$,  and $j(h)(\kappa)=\name{\mathbb{R}}$.
The iteration $j(\mathbb{P}_\kappa)$ is $\mathbb{P}'=\l \mathbb{P}'_\beta,\name{\mathbb{Q}}_\alpha: \beta\leq j(\kappa),\alpha<j(\kappa)\r$.
Since the critical point of $j$ is $\kappa$, $\mathbb{P}'\upharpoonright \kappa=\mathbb{P}_\kappa$, and since $j(h)(\kappa)=\name{\mathbb{R}}$, by definition of the iteration $\mathbb{P}'_{\kappa+1}=\mathbb{P}_\kappa\ast \name{\mathbb{R}}$.
In particular one can form the generic extension $M[G_{\mathbb{P}_\kappa}*G_{\mathbb{R}}]$ in the model $V[G_{\mathbb{P}_\kappa}*G_{\mathbb{R}}]$.

The tail of the iteration $\mathbb{P}'_{(\kappa+1,j(\kappa))}$ is now $\theta$-closed where $\theta$ is some $M[G_{\mathbb{P}_\kappa}*G_{\mathbb{R}}]$-inaccessible above $\kappa$.
Since $(2^\kappa)^{M}=2^\kappa$, it follows that $\mathbb{P}'_{(\kappa+1,j(\kappa))}$ is at least $(2^\kappa)^{+}$-closed from the point of view of $V$.
In $V$, consider the sets $\mathcal{C}=\{\name{C}:\name{C}\text{ is an } \mathbb{S}-\text{name for a club at }\kappa\}$ and $\mathcal{A}=\{\name{A}: \name{A}\text{ is an } \mathbb{S}-\text{name for a subset of }\kappa\}$.
By reducing to nice names we conclude that the cardinality of $\mathcal{C}$ and $\mathcal{A}$ is at most $2^{\kappa}$.
Since $M$ is closed under $2^\kappa$ sequences, it follows that $\{j(\name{C}):\name{C}\in \mathcal{C}\},\{j(\name{A}):\name{A}\in \mathcal{A}\}\in M$.

We claim that there is a condition $p\in \mathbb{P}'_{(\kappa+1,j(\kappa))}$ and $\delta$ such that $$p\Vdash_{\mathbb{P}'_{(\kappa+1,j(\kappa))}} \delta\in j(\name{S})\cap(\bigcap \{j(\name{C}):\name{C}\in \mathcal{C}\}).$$
Let us prove this claim.
Since $j(\kappa)$ is above $2^{\kappa}$ and $\{j(\name{C}):\name{C}\in \mathcal{C}\}$ is a collection $2^\kappa$-many names for clubs at $j(\kappa)$, it is forced by the empty condition that the intersection of these clubs is a club.
Since $j(\name{S})$ is a name for a stationary set in $j(\kappa)$, one can find an ordinal $\delta$ and a condition $p$ which forces that $\delta\in j(\name{S})\cap(\bigcap \{j(\name{C}):\name{C}\in \mathcal{C}\})$.
Now using the closure of the forcing we can find a single condition $p\leq p'\in \mathbb{P}'_{(\kappa+1,j(\kappa))}$ such that $p'$ decides the statement $\delta\in j(\name{A})$ for every $\name{A}\in \mathcal{A}$.

In $V[G_{\mathbb{P}_\kappa}\ast G_{\mathbb{R}}]$ define $$\mathscr{U}=\{\name{A}[G_{\mathbb{P}_\kappa}*G_{\mathbb{R}}]: \name{A}\in \mathcal{A}\wedge p'\Vdash_{\mathbb{P}'_{(\kappa+1,j(\kappa))}} \delta\in j(\name{A})\}.$$
We claim that $\mathscr{U}$ is a $\kappa$-complete ultrafilter over $\kappa$ which extends $\mathscr{D}\cup\{S\}$.
Indeed, ${\rm crit}(j)=\kappa$ and $j$ is an elementary embedding, so it is clear that $\mathscr{U}$ is a $\kappa$-complete ultrafilter over $\kappa$.
Since we have defined the condition $p'$ so that $p'\Vdash \delta\in j(\name{S})$ and $\delta\in j(\name{C})$ for every name of a club $\name{C}$, it follows by definition that $S=\name{S}[G_{\mathbb{P}_\kappa}*G_{\mathbb{R}}]\in \mathscr{U}$ and for every club $C\in V[G_{\mathbb{P}_\kappa}*G_{\mathbb{R}}]$, there is a name $\name{C}$ for $C$, and $p'\Vdash \delta\in j(\name{C})$.

\hfill \qedref{thmslimsupercompact}

Gathering everything under the same canopy, we can phrase the following:

\begin{corollary}
\label{cor2}
It is consistent to have a $\kappa$-complete ultrafilter over a measurable cardinal $\kappa$ which extends $\mathscr{D}_\kappa$ and fails to satisfy the Galvin property.
\end{corollary}

\hfill \qedref{cor2}

\newpage

\section{Open problems}

We conclude the paper with a short list of open problems related to our study.
Let $\mathscr{U}$ be a $\kappa$-complete filter over $\kappa$.
In \cite{bg}, the Galvin property was used in order to prove that quotients of Prikry-Magidor forcings with normal ultrafilters are $\kappa^+$-cc
in the generic extension, and therefore do not add fresh subsets to cardinals of cofinality greater than $\kappa$.
In this paper we proved the consistency of the existence of a $\kappa$-complete ultrafilter $\mathscr{U}$ which fails to satisfy the Galvin property, and hence the properties of the quotients of the tree Prikry forcing $\mathbb{P}_{\mathscr{T}}(\mathscr{U})$ are unclear:

\begin{question}
Is there a quotient of $\mathbb{P}_{\mathscr{T}}(\mathscr{U})$ which in not $\kappa^+$-cc in the full generic extension by $\mathbb{P}_{\mathscr{T}}(\mathscr{U})$?
\end{question}

One may ask, directly, about fresh sets:

\begin{question}
Is there a fresh subset of $\kappa^+$ after forcing with $\mathbb{P}_{\mathscr{T}}(\mathscr{U})$ with respect to some intermediate model?
\end{question}

As mentioned in the introduction, a product of $p$-points satisfies the Galvin property.
By the main result of this paper, it is consistent that there is a $\kappa$-complete ultrafilter over $\kappa$ which fails to satisfy the Galvin property.
Thus, such an ultrafilter cannot be a $p$-point or a product of $p$-points.
One may wonder whether there are more constraints on such ultrafilters.

\begin{question}
Are there more combinatorial properties of ultrafilters which imply the Galvin property?
\end{question}

The method of this paper fits very well the case of measurable cardinals, since such cardinals carry $\kappa$-complete ultrafilters.
It would be interesting to investigate a parallel problem at successor strongly regular cardinals.
The following is a modification of an open problem from \cite{bg}.

\begin{question}
\label{qgch} Suppose that $2^\kappa=\kappa^+$.
\begin{enumerate}
\item [$(\aleph)$] Is the statement $\neg{\rm Gal}(\mathscr{F},\kappa^+,\kappa^{++})$ consistent for some $\kappa^+$-complete filter $\mathscr{F}$ over $\kappa^+$?
\item [$(\beth)$] Is the statement $\neg{\rm Gal}(\mathscr{F},\kappa^+,\kappa^{++})$ consistent where $\mathscr{F}$ is $\kappa^+$-complete and $\mathscr{D}_{\kappa^+}\subseteq\mathscr{F}$?
\end{enumerate}
\end{question}

A natural model in which $\kappa$ is measurable and every $\kappa$-complete ultrafilter over $\kappa$ satisfies the Galvin property is $L[\mathscr{U}]$.
But the tolerance of inner models with respect to large cardinals is limited.

\begin{question}
\label{qsupercompact} Is it consistent that $\kappa$ is supercompact and ${\rm Gal}(\mathscr{U},\kappa,\kappa^+)$ holds at every $\kappa$-complete ultrafilter $\mathscr{U}$ over $\kappa$?
\end{question}

And a similar problem from a different perspective:

\begin{question}
\label{qglobal} Is it consistent with the existence of a measurable cardinal that every measurable cardinal carries a $\kappa$-complete ultrafilter which fails to satisfy the Galvin property?
In particular, is it consistent that there is a proper class of measurable cardinals and each one of them carries a $\kappa$-complete ultrafilter which fails to satisfy the Galvin property?
\end{question}

\section{Acknowledgements}

We thank the referee of the paper for many helpful suggestions and mathematical corrections.
We also thank F. Galvin for an inspiring work which stimulated the results presented here.

\newpage

\bibliographystyle{amsplain}
\bibliography{arlist}

\providecommand{\bysame}{\leavevmode\hbox to3em{\hrulefill}\thinspace}
\providecommand{\MR}{\relax\ifhmode\unskip\space\fi MR }
\providecommand{\MRhref}[2]{%
  \href{http://www.ams.org/mathscinet-getitem?mr=#1}{#2}
}
\providecommand{\href}[2]{#2}
\begin{thebibliography}{10}

\bibitem{MR830084}
U.~Abraham and S.~Shelah, \emph{On the intersection of closed unbounded sets},
  J. Symbolic Logic \textbf{51} (1986), no.~1, 180--189. \MR{830084}

\bibitem{MR0369081}
J.~E. Baumgartner, A.~H\c{a}j\c{n}al, and A.~Mate, \emph{Weak saturation
  properties of ideals}, Infinite and finite sets ({C}olloq., {K}eszthely,
  1973; dedicated to {P}. {E}rd\H{o}s on his 60th birthday), {V}ol. {I}, 1975,
  pp.~137--158. Colloq. Math. Soc. J\'{a}nos Bolyai, Vol. 10. \MR{0369081}

\bibitem{ezekiel}
Ezekiel Ben~Buzi, \emph{The book of {E}zekiel}, Prophets, 586 B.C.E.

\bibitem{bgp}
Tom Benhamou, Shimon Garti, and Alejandro Poveda, \emph{Negating the {G}alvin
  property at successors of singular cardinals}, preprint (2021).

\bibitem{bg}
Tom Benhamou and Moti Gitik, \emph{{I}ntermediate {M}odels of {M}agidor-{R}adin
  {F}orcing-{P}art {I}{I}}, Annals of Pure and Applied Logic \textbf{173}
  (2022), 103107.

\bibitem{MR2768691}
James Cummings, \emph{Iterated forcing and elementary embeddings}, Handbook of
  set theory. {V}ols. 1, 2, 3, Springer, Dordrecht, 2010, pp.~775--883.
  \MR{2768691}

\bibitem{MR750828}
Keith~J. Devlin, \emph{Constructibility}, Perspectives in Mathematical Logic,
  Springer-Verlag, Berlin, 1984. \MR{750828}

\bibitem{MR3604115}
Shimon Garti, \emph{Weak diamond and {G}alvin's property}, Period. Math.
  Hungar. \textbf{74} (2017), no.~1, 128--136. \MR{3604115}

\bibitem{MR3787522}
\bysame, \emph{Tiltan}, C. R. Math. Acad. Sci. Paris \textbf{356} (2018),
  no.~4, 351--359. \MR{3787522}

\bibitem{MR3902803}
Moti Gitik, \emph{Strange ultrafilters}, Arch. Math. Logic \textbf{58} (2019),
  no.~1-2, 35--52. \MR{3902803}

\bibitem{MR281605}
A.~Hajnal, \emph{On some combinatorial problems involving large cardinals},
  Fund. Math. \textbf{69} (1970), 39--53. \MR{281605}

\bibitem{MR277346}
Kenneth Kunen, \emph{Some applications of iterated ultrapowers in set theory},
  Ann. Math. Logic \textbf{1} (1970), 179--227. \MR{277346}

\bibitem{MR0472529}
Richard Laver, \emph{Making the supercompactness of {$\kappa $} indestructible
  under {$\kappa $}-directed closed forcing}, Israel J. Math. \textbf{29}
  (1978), no.~4, 385--388. \MR{0472529 (57 \#12226)}

\bibitem{MR2987148}
Philipp L\"{u}cke, \emph{{$\Sigma^1_1$}-definability at uncountable regular
  cardinals}, J. Symbolic Logic \textbf{77} (2012), no.~3, 1011--1046.
  \MR{2987148}

\bibitem{MR3075383}
Neil~H. Williams, \emph{Combinatorial set theory}, Studies in Logic and the
  Foundations of Mathematics, vol.~91, North-Holland Publishing Co., Amsterdam,
  1977. \MR{3075383}

\end{thebibliography}

\end{document}